# Admissible Region of Large-Scale Uncertain Wind Generation Considering Small-Signal Stability of Power Systems

Yanfei Pan, Shengwei Mei, *Fellow, IEEE,* Wei Wei, *Member, IEEE,* Chen Shen, *Senior Member, IEEE,* Jiabing Hu, *Senior Member, IEEE,* Feng Liu*, *Member, IEEE*

*Abstract*—The increasing integration of wind generation has brought great challenges to small-signal stability analysis of bulk power systems, since the volatility and uncertainty nature of wind generation may considerably affect equilibriums of the systems. In this regard, this paper develops a conceptual framework to depict the influence of uncertain wind power injections (WPIs) on small-signal stability of bulk power systems. To do this, a new concept, the *admissible region* of uncertain wind generation considering small-signal stability (*SSAR*) is introduced to geometrically measure how much uncertain wind generation can be accommodated by a bulk power system without breaking its small-signal stability. As a generalization of the traditional concept of the small-signal stability region (SSSR), SSAR is derived by extending the SSSR to a higher-dimensional injection space that incorporates both the conventional nodal generation injections and the WPIs, and then mapping it onto the lower-dimensional WPI space. Case studies on the modified New England 39-bus system with multiple wind farms illustrate the SSAR concept and its potential applications.

*Index Terms*—Power system, small-signal stability, wind generation, uncertainty, admissible region.

## I. Introduction

INCREASING integration of wind generation raises great concerns about power system stability, particularly the small-signal stability. Generally speaking, the new power source affects the small-signal stability of the bulk power system from two aspects: 1) the dynamic of wind turbine generators (WTGs) differ from that of the synchronous generators (SGs), which may produce new oscillation modes that interact with dynamics of the SGs; 2) the power supplied by WTGs is regarded as uncertain and non-dispatchable injection, as it is physically determined by volatile wind, and is accommodated in a "must-taken" manner. Uncertain variation of wind power injections (WPIs) may change the equilibriums of power systems, consequently influencing their small-signal stability. Practically, the uncertainty of WPIs can be approximately depicted by wind generation forecast error.

From the viewpoint of dynamic aspects, great efforts have been devoted to investigate how the dynamic of WTGs influence the small-signal stability. Ref. [1] characterizes the small-signal dynamic behaviors of DFIG as well as the effects of DFIG parameters by performing modal analysis on an SMIB system with a DFIG. In [2], the effects of DFIG on oscillation modes of a bulk power system are analyzed by replacing SGs with DFIGs. In [3], it is theoretically analyzed and numerically demonstrated that, under ideal conditions, the dynamic of a DFIG contributes little to the dominant modes of the original power system. Here, the term of "ideal conditions" means that the DFIG is controlled with the maximum power point tracking (MPPT) mode and the dynamic of phase lock loop (PLL) is neglected in DFIG model. In such circumstances, the influence of integration of wind farms on the dominant oscillations of the bulk power system is mainly boiled down to the equilibrium drift due to volatile and uncertain WPIs.

As for the wind generation uncertainty aspect, the Monte Carlo based method, the probability analysis based method, and the stochastic differential equation (SDE) based method are developed [4-6]. In [4] effects of wind power uncertainty on small-signal stability are analyzed by performing massive Monte Carlo simulations. The resulted distribution density of critical eigenvalues indicates the probabilistic stability of the power system. In [5], the probabilistic density function (PDF) of critical modes are derived directly from the PDF of multiple sources of wind generation, providing a systematic method to evaluate the influence of high-penetration wind generation on power system's small-signal stability. In [6], the mechanical power input of a wind turbine is regarded as a stochastic excitation to the system, leading to anSDE formulation of dynamic power system with uncertain wind generation disturbances. Then the standard SDE theory can be employed to investigate the impact of the stochastic excitation generated by wind generation on small-signal stability of power systems.

It is worthy of noting that, the aforementioned works are of point-wise fashion, where both the small-signal stability and the impact of wind generation are investigated in state space and rely on a given working point. This paper alternatively

This work was supported in part by the National Natural Science Foundation of China (No. 51377092, No. 51321005) and the National Basic Research Program (973 Program) of China (No. 2012CB215103)

Y. Pan, F. Liu, W. Wei, C. Shen and S. Mei are with the State Key Laboratory of Power Systems, Department of Electrical Engineering and Applied Electronic Technology, Tsinghua University, 100084 Beijing, China. (Corresponding author: Feng Liu, e-mail: lfeng@mail.tsinghua.edu.cn).

J. Hu is with the School of Electrical and Electronic Engineering and the State Key Laboratory of Advanced Electromagnetic Engineering and Technology, Huazhong University of Science and Technology, Wuhan 430074, China. (e-mail: j.hu@mail.hust.edu.cn).

investigates the problem from a region-wise fashion. To do that, we directly analyze the influences of wind generation from the perspective of small-signal stability region (SSSR).

Other than traditional eigen-analysis, the SSSR is defined on parameter space or power injection space, which depicts the feasibility region subjected to small-signal stability conditions. References [7] and [8] derive sufficient conditions for steady-state security regions to be small-signal stable. Later [9-11] point out that SSSR's boundary is composed by points of Hopf bifurcation (HB), saddle-node bifurcation (SNB) and singularity induced bifurcation (SIB), where HB is closely related to power system oscillations. [12-16] further propose efficient algorithms to compute SSSR boundaries efficiently. This paper aims to extend the existing concept of SSSR to cope with uncertain wind power injections. Particularly, two critical questions are considered: 1) how to find the limits (boundaries) of uncertain wind generation, within which the small-signal stability of system can be guaranteed; 2) how to quantitatively assess the impact of uncertain wind generation on small-signal stability from a region-wise point of view.

To answer these two questions, this paper first proposes a concept of admissible region of wind generation considering small-signal stability (SSAR) that is defined on the space of wind power injections (WPIs). It mathematically depict the region within which uncertain WPIs varies without breaking the small-signal stability of the bulk power system. Then the uncertainty of wind generation is modelled as an ellipsoidal uncertainty set (EUS). By checking whether or not the uncertainty set is completely inside the inner of the SSAR, we can directly judge if there is certain probability that the WPIs may cause small-signal instability, and how much the probability is. This essentially provides a geometric description for the capability of power system to accommodate uncertain wind generation, enabling quantitative assessment of the small-signal stability under wind generation uncertainty in an intuitive and visual fashion.

The remainder of this paper is organized as follows. Section II introduces the SSSR theory and generalizes it to incorporate wind generation. The concept of SSAR is proposed in Section III. Section IV addresses the modelling of wind generation uncertainty based on the EUS. Then potential applications of the SSAR are discussed in Section V. In Section VI, a modified New England 39-bus system with three wind farms is used to perform case studies. Finally, Section VII concludes the paper.

## II. GENERALIZING SMALL-SIGNAL STABILITY REGION

### A. Preliminary of small-signal stability region

Mathematically, a power system can be described by a set of differential-algebraic equations (DAEs) with parameters:

$$\begin{cases} \dot{x} = F(x, y, p) \\ 0 = G(x, y, p) \end{cases} \quad (1)$$

where $x \in \mathbb{R}^n$, $y \in \mathbb{R}^m$ are vectors of state and algebraic variables, respectively. $p \in \mathbb{R}^l$ is the vector of parameters.

Mathematically, (1) can be linearized at a given equilibrium point, yielding the following augmented state equation:

$$\begin{bmatrix} \Delta \dot{x} \\ 0 \end{bmatrix} = \begin{bmatrix} \tilde{A}(p) & \tilde{B}(p) \\ \tilde{C}(p) & \tilde{D}(p) \end{bmatrix} \begin{bmatrix} \Delta x \\ \Delta y \end{bmatrix} \quad (2)$$

Without loss of clarity, the parameter "$p$" is omitted for simplification in the following parts. In (2), $\tilde{A} = \partial F / \partial x \in \mathbb{R}^{n \times n}$, $\tilde{B} = \partial F / \partial y \in \mathbb{R}^{n \times m}$, $\tilde{C} = \partial G / \partial x \in \mathbb{R}^{m \times n}$ and $\tilde{D} = \partial G / \partial y \in \mathbb{R}^{m \times m}$. Assume that $\tilde{D}$ is nonsingular, then (2) can be reduced into the following form

$$\Delta \dot{x} = A \Delta x \quad (3)$$

where $A = \tilde{A} - \tilde{B}\tilde{D}^{-1}\tilde{C}$ is the reduced state matrix. $\lambda = [\lambda_1, \lambda_2, \cdots, \lambda_{2n}]$ is the eigenvalue vector of $A$.

It is well known that, when each $\lambda_i \in \lambda$ have negative real parts, the system is small-signal stable. Note that, both $A$ and $\lambda$ are parameterized by $p$. Then the parameter vector $p$ can span an $l$-dimensional space, $S_p := \text{span}\{p_1, p_2, \cdots, p_l\}$, which is referred to as a parameter space. Consequently, the SSSR, denoted by $\Omega_{\text{SSSR}}$, can be defined on the parameter space [15]:

$$\Omega_{\text{SSSR}} := \left\{ p \, \middle| \, \begin{array}{l} \max\{\text{Re}(\lambda_i)\} < 0 \ \forall \lambda_i \in \lambda, \ \forall i \\ \tilde{D} \text{ is nonsingular,} \end{array} \right\} \quad (4)$$

where $\text{Re}(\lambda_i) = \text{real}(\lambda_i)$ is the real part of $\lambda_i$. Then the boundary of SSSR, denoted by $\partial \Omega_{\text{SSSR}}$, can be further defined as:

$$\partial \Omega_{\text{SSSR}} := \left\{ p \, \middle| \, \begin{array}{l} \max\{\text{Re}(\lambda_i)\} = 0 \\ \forall \lambda_i \in \lambda, \forall i \end{array} \right\} \cup \left\{ p \, \middle| \, \tilde{D} \text{ is singular} \right\} \quad (5)$$

where $\max\{\text{Re}(\lambda_i)\} = 0$ refers to Hopf bifurcation (HB) or saddle-node bifurcation (SNB), and the singularity of $\tilde{D}$ refers to singularity induced bifurcation (SIB).

Usually, $p$ can be selected as the active power injection vector of $n_s$ generator nodes $p_s = [p_{s1}, p_{s2}, \cdots, p_{sn_s}]^T$. Assuming that the power loss of the system is neglected and the system load is constant as $L$, to keep the power balance, the power injections of SGs must satisfy

$$p_{s1} + p_{s2} + \cdots + p_{sn_s} = L \quad (6)$$

Then the corresponding parameter space is $S_{p_s} := \text{span}\{p_{s1}, p_{s2}, \cdots, p_{sn_s}\}$, which is referred to as SGs' power injection space. The SSSR defined on this power injection space is denoted by $\Omega_{\text{SSSR}\_p_s}$. It can be calculated off-line and utilized by system operators on-line. The relative position of an operating point to the SSSR boundaries as well as the distance between them can provide useful information for operators to make decisions in operation.

### B. Extending SSSR to incorporate WPIs

As well known, the integration of wind power generation creates great challenges to the small-signal stability analysis. The SSSR analysis is in more difficult case as it is concerned with impact of uncertain wind generation on the entire small-signal stability region as well as its boundaries, other than a certain specific working point. Thus the first step is to extend the conventional SSSR concept so as to enable the consideration of uncertain wind power generation. As mentioned in the introduction, we concentrate on the influences of the uncertain WPIs while neglecting the dynamics of WTGs. It is worth noting that, this simplification is not necessary because the proposed methodology is generic and does not rely on concrete dynamic system model.

*1) Generator regulation to cope with WPIs:* Before extending the SSSR concept, we need to address how the power balance of power system is attained when volatile wind generation is injected. Physically, whenever variation of WPIs causes power unbalance, the SGs or other controllable power sources will be adjusted either manually or automatically to cope with WPIs in real time. It follows

$$\sum_{i=1}^{n_s}(p_{si0}+\Delta p_{si})+\sum_{j=1}^{n_w}(p_{wj0}+\Delta p_{wj})=L$$
$$p_{si}^l \leq p_{si0}+\Delta p_{si} \leq p_{si}^u \quad (7)$$
$$p_{wj}^l \leq p_{wj0}+\Delta p_{wj} \leq p_{wj}^u$$

where $n_w$ is the number of WPIs. The superscript "l" and "u" denote the lower and upper limit, respectively. The subscript "0" means the forecast or scheduled power injection. $\Delta p_{wj}$ refers to the uncertain part of WPI$_j$, while $\Delta p_{si}$ is the regulation of the *i*th nodal power injection of SGs to deal with the power unbalance caused by $\Delta p_{wj}$. This can be achieved by deploying automatic generation control (AGC) in operation. For simplicity we ignore the dynamics of the AGC, and the mismatched power is directly distributed to SGs with a given contribution factor vector $\boldsymbol{\gamma}=(\gamma_1,\gamma_2,\cdots,\gamma_{n_s})$ satisfying $\sum_{i=1}^{n_s}\gamma_i=1$. Thus, we have

$$\Delta p_{si}=\gamma_i\sum_{j=1}^{n_w}\Delta p_{wj} \quad (8)$$

*2) Extending SSSR to incorporate WPIs:* To incorporate non-dispatchable WPIs, the SGs' power injection space, $\mathcal{S}_{\mathbf{p}_s}$, needs to be expanded into a higher-dimensional nodal injection space, denoted by $\mathcal{S}_{\mathbf{p}_e}:=\mathrm{span}\{p_{s1},\cdots,p_{sn_s},p_{w1},\cdots,p_{wn_w}\}$. It is spanned by both conventional SGs' power injections and WPIs. Then the extended SSSR can be defined on $\mathcal{S}_{\mathbf{p}_e}$, which is

$$\Omega_{\mathrm{SSSR\_p_e}}:=\left\{\boldsymbol{p}_e \left| \begin{array}{l} \max\{\mathrm{Re}(\lambda_i)\}<0 \ \forall \lambda_i \in \boldsymbol{\lambda}, \ \forall i \\ \tilde{\boldsymbol{D}} \text{ is nonsingular} \\ \text{constraint (7)} \end{array} \right. \right\} \quad (9)$$

where $\boldsymbol{p}_e=[p_{s1},\cdots,p_{sn_s},p_{w1},\cdots,p_{wn_w}]^\mathrm{T}$.

To further take into account the effect of AGC, constraint (8) should also be augmented. Note that the AGC regulation is always associated with a given contribution factor vector $\boldsymbol{\gamma}$. Denote the $\Omega_{\mathrm{SSSR\_p_e}}$ under a certain $\boldsymbol{\gamma}$ by $\Omega_{\mathrm{SSSR\_p_e}}^{\boldsymbol{\gamma}}$, which is referred to as $\boldsymbol{\gamma}$-SSSR. Then it can be defined as

$$\Omega_{\mathrm{SSSR\_p_e}}^{\boldsymbol{\gamma}}=\left\{\boldsymbol{p}_e \left| \begin{array}{l} \max\{\mathrm{Re}(\lambda_i)\}<0 \ \forall \lambda_i \in \boldsymbol{\lambda}, \ \forall i \\ \tilde{\boldsymbol{D}} \text{ is nonsingular;} \\ \text{(7) and (8) for a given } \boldsymbol{\gamma} \end{array} \right. \right\} \quad (10)$$

Obviously, there is $\Omega_{\mathrm{SSSR\_p_e}}^{\boldsymbol{\gamma}} \subset \Omega_{\mathrm{SSSR\_p_e}}$, which seems to introduce conservativeness. However, since AGC mode is always determined ex ante, this treatment is reasonable and makes sense for system operators.

Note that $\Omega_{\mathrm{SSSR\_p_e}}^{\boldsymbol{\gamma}}$ can be regarded as a special case of the SSSR on power injection space complying with operation rules and constraints. Since the constraints are linear, its main characteristics can inherit from that of SSSR. In [17], some topology characteristics of SSSR are discussed, such as the existence of holes (instability region) inside SSSR. Fortunately, it is also revealed that the hole occurs only due to degeneration of Hopf bifurcations, which is mainly resulted from inappropriate exciter parameters. To facilitate the subsequent analysis we make the assumption that all exciter parameters are in appropriate ranges such that there is no holes inside $\Omega_{\mathrm{SSSR\_p_e}}^{\boldsymbol{\gamma}}$.

III. ADMISSIBLE REGION OF WIND GENERATION CONSIDERING SMALL-SIGNAL STABILITY

*A. Definition*

In a bulk power system with WPIs, it is crucial to mathematically depict the amount of wind power that can be accommodated by the system, without breaking the small-signal stability conditions. To do this, a metric defined in the space of WPIs is desirable. In this regard, we introduce a new concept, *admissible region* of wind generation considering small-signal stability (SSAR), which is defined on the WPI space. Let $\boldsymbol{p}_e=[p_{s1},\cdots,p_{sn_s},p_{w1},\cdots,p_{wn_w}]^\mathrm{T}$ be the vector of nodal power injections comprising both the conventional SGs' power injections and the WPIs. Then the SSAR, denoted by $\Omega_{\mathrm{SSAR\_p_w}}$ is defined as follows.

**Definition 1**: SSAR is the region defined on the WPI space, which satisfies

$$\Omega_{\mathrm{SSAR\_p_w}}:=\left\{\boldsymbol{p}_w \in \mathbb{R}^{n_w} \left| \exists \boldsymbol{p}_s \text{ such that } \Omega_{\mathrm{SSAR\_p_e}} \neq \varnothing \right. \right\} \quad (11)$$

where, $\boldsymbol{p}_s=[\boldsymbol{I}_{n_s} \ \boldsymbol{0}]\boldsymbol{p}_e$ represents the vector of conventional SGs' power injections, while $\boldsymbol{p}_w=[\boldsymbol{0} \ \boldsymbol{I}_{n_w}]\boldsymbol{p}_e$ is the vector of WPIs.

Furthermore, when considering AGC regulation with distribution factor $\boldsymbol{\gamma}$, the SSAR, denoted by $\Omega_{\mathrm{SSAR\_p_w}}^{\boldsymbol{\gamma}}$ can be defined as:

**Definition 2**: SSAR with an AGC distribution factor $\boldsymbol{\gamma}$ ($\boldsymbol{\gamma}$-SSAR) is the region defined on the WPI space satisfying

$$\Omega_{\mathrm{SSAR\_p_w}}^{\boldsymbol{\gamma}}:=\left\{\boldsymbol{p}_w \in \mathbb{R}^{n_w} \left| \boldsymbol{p}_e \in \Omega_{\mathrm{SSSR\_p_e}}^{\boldsymbol{\gamma}} \right. \right\} \quad (12)$$

This definition indicates that a WPI vector $\boldsymbol{p}_w$ is small-signal stable, if the conventional power injections under the given $\boldsymbol{\gamma}$ can retain the operating point still inside $\Omega_{\mathrm{SSSR\_p_e}}^{\boldsymbol{\gamma}}$. However, as mentioned previously, it is more practical to use $\boldsymbol{\gamma}$-SSAR than SSAR since AGC regulation is always required in power system operation. Thus, this paper will focus on the $\boldsymbol{\gamma}$-SSAR.

*B. Constructing SSAR based on extended SSSR*

From the definitions of SSSR and SSAR, it can be found that SSAR is the projection of SSSR to the lower-dimensional WPI space. That is:

$$\Omega_{\mathrm{SSAR\_p_w}}=\mathrm{proj}_{\boldsymbol{p}_w}(\Omega_{\mathrm{SSSR\_p_e}}) \quad (13)$$

where proj(·) is the projection operator. Similarly, we have

$$\Omega_{\mathrm{SSAR\_p_w}}^{\boldsymbol{\gamma}}=\mathrm{proj}_{\boldsymbol{p}_w}(\Omega_{\mathrm{SSSR\_p_e}}^{\boldsymbol{\gamma}}) \quad (14)$$

Note that $\Delta p_{si}$ can be obtained from $\Delta p_{wj}$ according to (8). Then SSAR and **γ**-SSAR can be directly obtained by eliminating $\Delta p_{si}$ in SSSR and **γ**-SSSR, respectively.

To illustrate this clearly, a simple system with one synchronous generator and two wind farms is taken as an example. $S_{p_e}:=\text{span}\{p_s,p_{w1},p_{w2}\}$ constitutes the power injection space. Assume that load $L$ is constant and the unbalanced power is eliminated by the synchronous generator, satisfying

$$p_{w10} + \Delta p_{w1} + p_{w20} + \Delta p_{w2} + p_{s0} + \Delta p_s = L \quad (15)$$

(15) indicates that the $\Omega_{SSSR\_p_e}$ is a region on a plane on the 3-D power injection space, $S_{p_e}$, determined by (4).

Fig.1 depicts $\Omega_{SSSR\_p_e}$ and $\Omega_{SSAR\_p_w}$. Operating point A in $\Omega_{SSAR\_p_w}$ is small-signal stable if there exists $p_s$ to make sure that the operating point is in $\Omega_{SSSR\_p_e}$, shown as point B.

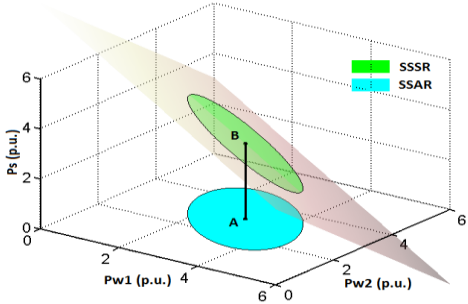

Fig. 1. $\Omega_{SSSR\_p_e}$ and $\Omega_{SSAR\_p_w}$

To take into account effect of AGC regulation, another system with two synchronous generators and two wind farms is also used for illustration. Here $S_{p_e}:=\text{span}\{p_{s1},p_{s2},p_{w1},p_{w2}\}$, the unbalanced power is distributed to two generators with a distribution factor vector, **γ**=($\gamma_1,\gamma_2$), satisfying

$$\sum_{i=1}^{2}(p_{wi0}+\Delta p_{wi}) + \sum_{j=1}^{2}(p_{sj0}+\Delta p_{sj}) = L$$
$$\Delta p_{s1} = -\gamma_1(\Delta p_{w1}+\Delta p_{w2}) \quad (16)$$
$$\Delta p_{s2} = -\gamma_2(\Delta p_{w1}+\Delta p_{w2})$$

Fig.2 depicts two different $\Omega^{\gamma}_{SSAR\_p_w}$ with $\gamma_1$ and $\gamma_2$, on the plane spanned by $p_{w1}$ and $p_{w2}$. It can be observed that $\Omega^{\gamma}_{SSAR\_p_w}$ with $\gamma_2$ is larger than that with $\gamma_1$. This implies that an optimal **γ** may enhance power system's small-signal stability.

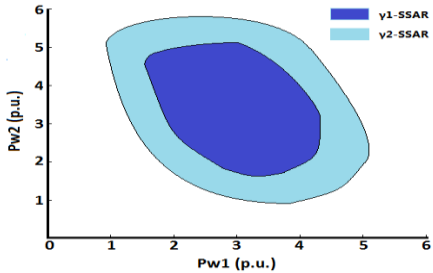

Fig. 2. $\Omega^{\gamma}_{SSAR\_p_w}$ corresponding to $\gamma_1$ and $\gamma_2$

### C. Computing SSAR boundaries

Mathematically, SSAR boundary is described by a set of implicit nonlinear equations, and it is difficult to obtain the analytical expression. In [16], a polynomial approximation method is proposed for computing SSSR boundaries based on the implicit function theorem. This can be directly applied to compute boundaries of the extended SSSR. Specifically, in terms of **γ**-SSSR, the quadratic approximation of $\partial\Omega^{\gamma}_{SSSR\_p_e}$ associated with the expansion point $p_{e0}$ can be expressed as

$$\Delta p_{e1} = \sum_{i=2}^{n_w+n_s}\frac{\partial p_{e1}}{\partial p_{ei}}\Delta p_{ei} + \frac{1}{2}\sum_{i=2}^{n_w+n_s}\sum_{j=2}^{n_w+n_s}\frac{\partial^2 p_{e1}}{\partial p_{ei}\partial p_{ej}}\Delta p_{ei}\Delta p_{ej} \quad (17)$$

where $\Delta p_{ei}=p_{ei}-p_{ei0}$, ($i=1,2,\cdots,n_s+n_w$) compromises both SGs' nodal power injections ($i\leq n_s$) and WPIs ($n_s<i\leq n_s+n_w$).

The approximated boundaries of **γ**-SSAR can be obtained by further eliminating $\Delta p_{ei}$ ($i\leq n_s$) in (17), which gives

$$\Delta p_{w1} = \sum_{i=2}^{n_w}\frac{\partial p_{w1}}{\partial p_{wi}}\Delta p_{wi} + \frac{1}{2}\sum_{i=2}^{n_w}\sum_{j=2}^{n_w}\frac{\partial^2 p_{w1}}{\partial p_{wi}\partial p_{wj}}\Delta p_{wi}\Delta p_{wj} \quad (18)$$

One can also obtain the approximate expression of the boundary by using fitting approach. It may provide the approximation of the boundary in a broader scope. However, enough number of critical points on the boundary are required.

Generally, the boundary of SSAR is difficult to visualize when $n_w$ is greater than three. Thus, in practice, SSAR is usually reduced to a lower-dimensional space (usually up to three) by fixing certain WPIs at their nominal values (values at the expansion point of the approximate boundary). This produces a series of profiles of the SSAR in different nodal injection subspaces.

### IV. MODELLING UNCERTAINTY OF WIND GENERATION

According to the concept of SSAR derived in the previous section, a forecast point of wind generation is admissible if it is located inside $\Omega^{\gamma}_{SSAR\_p_w}$. However, due to existence of forecast error, the realization of wind generation may deviate away from its forecast value. Thus, we need to carefully examine whether or not the deviation could cause the point turn to be out of $\Omega^{\gamma}_{SSAR\_p_w}$. As shown in Fig.3, the blue dots represent the distribution of forecast error of wind generation. It is observed that some of the dots are outside the SSAR, indicating certain possibility of small-signal instability of the power system. This simple example demonstrates that the SSAR can serve as a geometric measure to assess the level of small-signal stability of the power system with uncertain wind generation. To do this, we need to properly depict wind generation uncertainty.

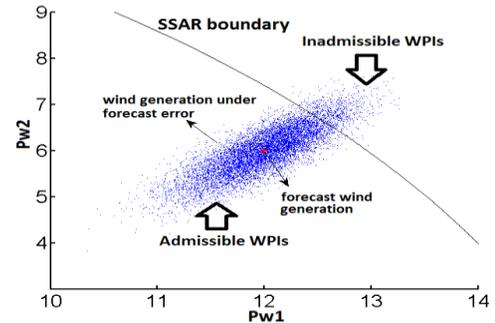

Fig. 3. Forecast point of wind generation inside SSAR and its possible deviations under certain forecast error

From the aforementioned example, the uncertainty of wind generation is modeled based on forecast error. This approach has been extensively investigated in unit commitment and economical dispatch. The uncertainty set, expressed by a set of inequalities [18-21], is a region covering certain confidence level (e.g., 90%) of possible realizations of a predict point under forecast error [21]. There are usually polyhedral sets [18,19] or ellipsoidal sets [21]. This paper uses the ellipsoidal uncertainty set (EUS) due to its concise expression and the ability to model correlations among different WPIs. It should be noted that our approach is general and other kinds of uncertainty sets also can apply.

The uncertain wind generation can be described by the following generalized ellipsoidal uncertainty set:

$$W^{\text{EUS}} = \left\{ p_{\mathbf{w}} \mid (p_{\mathbf{w}} - p_{\mathbf{w0}})^{\text{T}} Q^{-1} (p_{\mathbf{w}} - p_{\mathbf{w0}}) \le \eta \right\} \quad (19)$$

where, $p_{\mathbf{w}}$ is the WPI vector under forecast errors, $p_{\mathbf{w0}}$ the vector of forecast point. One can adjust the robustness of $W^{\text{EUS}}$ by altering the value of $\eta$. Denote by $\text{WFE}_i$ the wind forecast error of the $i$th WPI. Then the matrix $Q$ has the form of the covariance matrix:

$$Q = \begin{bmatrix} \sigma_1^2 & \rho_{12}\sigma_1\sigma_2 & \cdots & \rho_{1n_w}\sigma_1\sigma_{n_w} \\ \rho_{21}\sigma_1\sigma_2 & \sigma_2^2 & \cdots & \rho_{2n_w}\sigma_2\sigma_{n_w} \\ \vdots & \vdots & \ddots & \vdots \\ \rho_{n_w 1}\sigma_1\sigma_{n_w} & \rho_{n_w 2}\sigma_2\sigma_{n_w} & \cdots & \sigma_{n_w}^2 \end{bmatrix} \quad (20)$$

where $\sigma_i$ is the standard deviation of $\text{WFE}_i$; $\rho_{ij}$ the correlation between $\text{WFE}_i$ and $\text{WFE}_j$. Note that the value of $\rho_{ij}$ in this paper is not the well-known linear correlation, but the rank correlation. The comparison between the two correlations and the advantages of the rank correlation are given in Section III.B of reference [22]. $\eta$ represents the uncertainty measure determined by the confidence level of preference. The larger $\eta$ is, the larger the $W^{\text{EUS}}$. Nevertheless, a too large $\eta$ may bring over-conservativeness to the admissibility assessment of uncertain wind generation. To minimize the conservativeness, the following optimization can be deployed

$$\min \eta$$
$$\text{s.t.} \quad \Pr\left\{ p_{\mathbf{w}}^{\text{sample}} \mid p_{\mathbf{w}}^{\text{sample}} \in P_{\mathbf{W}}^{\text{EUS}} \right\} \ge \alpha \quad (21)$$
$$W^{\text{EUS}} = \left\{ p_{\mathbf{w}} \mid (p_{\mathbf{w}} - p_{\mathbf{w0}})^{\text{T}} Q^{-1} (p_{\mathbf{w}} - p_{\mathbf{w0}}) \le \eta \right\}$$

where $p_{\mathbf{w}}^{\text{sample}}$ is the collection of samples under certain forecast errors, which can be historical data or scenarios generated following certain distribution function. "Pr" means probability. $\alpha$ is the confidence probability. (21) can be elucidated as minimizing $\eta$ while guaranteeing that the probability of $p_{\mathbf{w}}^{\text{sample}}$ located inside $W^{\text{EUS}}$ is not less than $\alpha$.

## V. POTENTIAL APPLICATIONS

The proposed SSAR provides a geometrical metric on the small-signal stability of the power system under uncertain wind generation, characterizing exactly how much uncertainty the system can accommodate without causing instability. It is useful in power system security assessment issues. Note that the SSAR can be computed off-line and applied on-line with no need of updates until the topology changes. The feature is appealing in on-line or real-time applications. Several potential applications are suggested.

*1) Small-signal stability assessment:* With the ellipsoidal uncertainty set presented above, it is easy to justify that the subset $W_s = W^{\text{EUS}} \cap \Omega_{\text{SSAR\_P}_w}^\gamma$ is admissible and ensures the small-signal stability of the power system, while the subset $W_u = W^{\text{EUS}} \setminus W_s$ is inadmissible as any point contained in $W_u$ corresponds to a small-signal unstable state of the power system. Thus, the SSAR provides a quantitative way to assess not only whether or not the bulk power system is small-signal stable with forecast wind generation, but also with an uncertain set of forecast error. The correlation among multiple WPIs can also be taken into account. Besides, if the uncertain wind generation cannot be fully accommodated, it can quantitatively measure the probability of the bulk power system being stable or unstable under the uncertainty.

*2) Stability margin and vulnerable direction of wind generation variation*: The EUS shrinks with decreasing confidence probability $\alpha$. In case the EUS intersects with the SSAR boundary, if one gradually reduces $\alpha$ until the boundary of EUS is tangent with the SSAR boundary, then the direction from the forecast point toward the tangent point, denoted by $\vec{D}_{\text{vulner}}$, can be regarded as the most vulnerable direction of wind generation variation. Accordingly, the distance between the forecast point and the tangent point provides the least margins of WPIs, which can serve as an indicator to perceive potential risky variations of wind generation.

*3) Robust small-signal stability region*: As the situation where the EUS is tangent with the SSAR boundary indicates a critical state that the system is small-signal stable under uncertainty, the binding forecast point can be seen as a boundary point of a region within which the system can withstand the uncertainty of wind generation and maintain stability. This new region is essentially a robust small-signal stability region against all uncertainty of wind generation under consideration. In system monitoring, if the forecast operating point is outside this region, the system has certain probability of losing stability, and preventive actions could be taken for the sake of enhancing operation security.

*4) Wind generation curtailment*: The part of the EUS outside the SSAR boundary could be eliminated by adjusting the wind generation, e.g. curtailment. However, it is a challenging task for system operators as there are numerous adjustable directions. The dimension-reduced EUS and SSAR boundaries on different nodal injection subspaces may facilitate finding out the feasible adjustment of wind generation to improve stability in a visualized fashion. This problem also can be formulated as an optimization problem to minimizing the amount of wind generation curtailment such that the stability requirement can be satisfied.

## VI. ILLUSTRATIVE EXAMPLES

In this section, the New England 39-bus system [23] is modified and employed to illustrate the concept of SSAR. Three wind farms $W_1$, $W_2$ and $W_3$ are connected to the grid at buses 34, 35, and 37, respectively. The capacities of $W_1$, $W_2$ and $W_3$ are 1500MW, 1000MW and 800MW, respectively. The base value is chosen as $S_B$=100MW. The wind farms are modelled by power injections. G2 (the swing bus) and G10 are modelled by the classic model, while G1 and G3-G9 are modelled by the 3rd-order model with the simplified 3rd-order exciter model demonstrated in Appendix A.

### A. Base case

The AGC distribution factor vector is chosen as $\gamma=[0,0,0.1,0.1,0.05,0.05,0.1,0.1,0.1,0.4]^T$ according to the SGs' capacities for regulation. $WFE_1$-$WFE_3$ are assumed to follow the beta distribution as many studies did [24-25]. The analytical expression of beta distribution is shown in Appendix B. The forecast wind generation for $W_1 \sim W_3$ are $(p_{w10}, p_{w20}, p_{w30}) = (12, 6, 2)$ p.u. with standard deviations $\sigma_1=0.4, \sigma_2=0.5, \sigma_3=0.4$, respectively. The rank correlation matrix of $WFE_1$-$WFE_3$ is given in Tab.I.

TABLE I
RANK CORRELATION MATRIX OF $WFE_1$-$WFE_3$

|     | W1  | W2  | W3  |
| --- | --- | --- | --- |
| W1  | 1   | 0.8 | 0.5 |
| W2  | 0.8 | 1   | 0.8 |
| W3  | 0.5 | 0.8 | 1   |

In this case, 10000 random points are generated with the forecast error distribution, in a similar way as presented in [25]. The EUS is generated by using (21), with $\eta=7.30$ determined by $\alpha=95\%$. The EUS and part of the $\gamma$-SSAR boundary are depicted in Fig.4, in which the red points denote the possible realization of wind generation (scenarios) under forecast error, and the green ellipsoid is the EUS. As mentioned in Part C of Section III, the approximate expression of the $\gamma$-SSAR boundary in 3-D space can be obtained by using polynomial approximation or surface fitting approaches. However, for the sake of obtaining the exact boundary, we simply adopt the searching algorithm here.

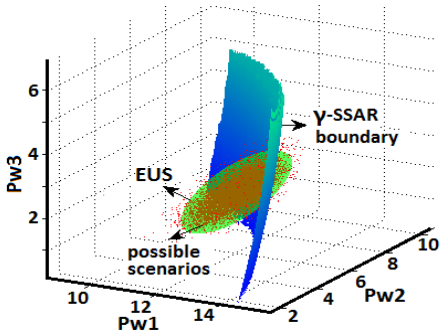

Fig. 4. The $\gamma$-SSAR boundary and the EUS of $WPI_{1-3}$

Fig.4 shows that some possible scenarios are located outside the $\gamma$-SSAR, and the EUS intersects with the $\gamma$-SSAR boundary. The instability probability under the uncertainty of wind generation, denoted by $P_{instab}$, could be calculated via integration approach as the expressions of the EUS and $\gamma$-SSAR boundary are both known. However, it is not easy as the expression of SSAR boundaries can be very complex. Therefore we alternatively adopt a Monte Carlo method to evaluate the instability probability using the equation below

$$P_{instab} = N_{out} / N_{total} \qquad (22)$$

where, $P_{instab}$ is the probability of small-signal instability; $N_{total}$ the number of total samples, $N_{out}$ the samples outside the $\gamma$-SSAR. In this case $N_{out}=666$, $N_{total}=10000$, thus $P_{instab}=6.66\%$. It should be emphasized that $N_{out}$ is obtained with no need of processing any eigen-analysis, but directly by substituting them into the expression of the SSAR boundary and then merely examining the sign of results.

### B. Effect of AGC distribution factor $\gamma$

As mentioned previously, the $\gamma$-SSAR boundary varies according to different $\gamma$. To illustrate this, the $\gamma$-SSAR under $\gamma_2=[0.1,0.1,0.1,0.1,0.05,0.05,0.1,0.1,0.3,0]^T$ is computed. This new distribution relieves G10 from the obligation of power regulation. The $\gamma_2$-SSAR boundary is then depicted and compared with that in the base case in Fig.5. The distribution factor vector in the base case is denoted as $\gamma_1$. As shown in Fig.5, $\gamma_2$-SSAR can accommodate more wind generation, while guaranteeing the small-signal stability of the system.

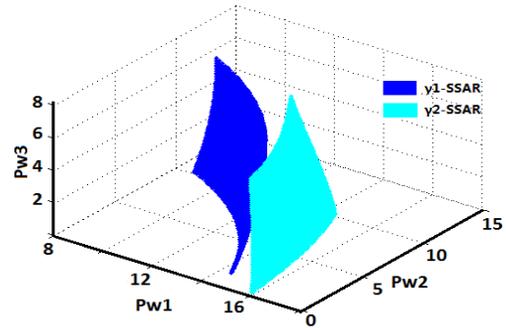

Fig. 5. the $\gamma$-SSAR boundary with $\gamma_1$ and $\gamma_2$

TABLE II
RANK CORRELATION MATRIX IN WEAKLY CORRELATED CASE

|     | W1  | W2  | W3  |
| --- | --- | --- | --- |
| W1  | 1   | 0.1 | 0.1 |
| W2  | 0.1 | 1   | 0.1 |
| W3  | 0.1 | 0.1 | 1   |

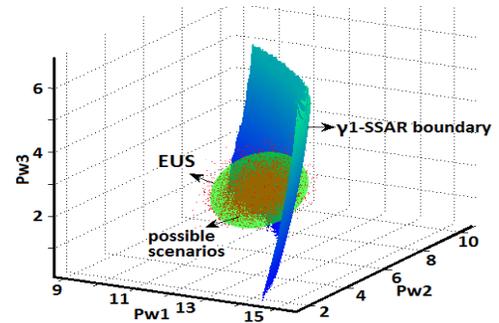

Fig. 6. the $\gamma_1$-SSAR boundary and the EUS of weakly correlated $WPI_{1-3}$

### C. Effect of WFE correlation

In the base case, $WFE_1$-$WFE_3$ are strongly correlated. Here, the case of weak correlation is studied. Assume that the correlation coefficient matrix is replaced by Tab. II, and other characteristics of the forecast error distribution are unchanged. 10000 random points and the EUS are generated in the same way as in the base case, shown with the $\gamma_1$-SSAR boundary in Fig.6. The resulted $P_{instab}$ is reduced to 3.71%. It should be pointed out that this does not implies that weak correlation of WFE must lead to lower $P_{instab}$. Our methodology, however, provides a quantitative way to assess the effect of correlation.

*D. Security assessments*

The security assessment is performed on 2-D spaces that is easy to visualize and most likely to apply in practical system monitoring. Assume that the correlations among different WFEs are the same as in the base case.

*1) Vulnerable direction of wind generation variation*: The EUS is tangent with the SSAR boundary when $\alpha$ decreases to 71.8%, as shown in Fig.7, where O is the forecast point, A is the point of tangency and B locates outside the level set with $\alpha$=71.8% but inside the SSAR. Obviously, point A has a smaller $\alpha$ than that of the point B, which means the system has a larger probability to operate at A than B. Since A is on the SSAR boundary which indicates a critical stability state while B is stable still with certain stability margin, it is reasonable to take direction $\overrightarrow{OA}$ as the most vulnerable direction of the possible variation of wind generation.

*2) Stability margin*: In Fig.7, it is easy to observe that the least margin of $W_1$ is 0.66 p.u. while that of $W_2$ is 0.55 p.u.. In practice, multiple profiles in different 2-D planes can be generated according to different WPIs of interest. Then the SSAR and its boundary are capable of providing system operators a visual tool to monitor the small-signal stability of the power system under uncertain wind generation.

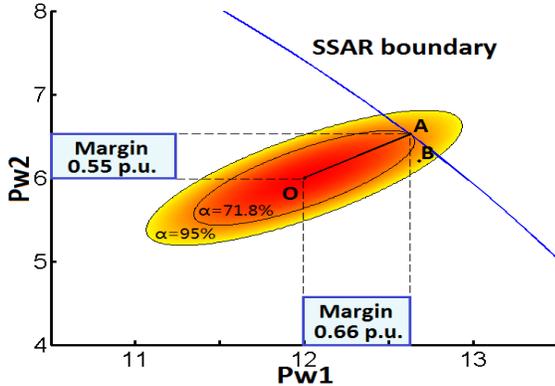

Fig. 7. The vulnerable direction and stability margin in different directions

*3) Wind generation adjustments*: When the instability probability, $P_{instab}$, for certain forecast point is too large to be acceptable, it is necessary to adjust the outputs of wind farms, say, wind spillage. The dimension-reduced SSARs and EUSs on different 2-D planes then are helpful to identify which WPIs are responsible to be adjusted to reduce $P_{instab}$ such that it turns to be acceptable. Fig.8 shows the dimension-reduced EUSs and the SSAR boundaries on three profiles in 2-D space: $(p_{w1},p_{w2})$, $(p_{w1},p_{w3})$, $(p_{w2},p_{w3})$. It is observed from Fig.8 that the EUS intersects with the SSAR boundary on the subspace of $(p_{w1},p_{w2})$, while in other two subspaces it stays inside the SSAR and has a certain stability margin. Thus, to reduce the probability of losing stability, $p_{w1}$ and $p_{w2}$ should be decreased. Furthermore, according to Fig.7, we speculate cautiously that decreasing $p_{w1}$ will improve the stability level more efficiently than decreasing $p_{w2}$, since the angle between the major axis of the EUS and the $p_{w1}$ axis is smaller. Meanwhile, we conjecture from Fig.8(b) that decreasing $p_{w3}$ would lower the stability level other than to enhance it as one's common sense, since it will move the EUS closer to the SSAR boundary.

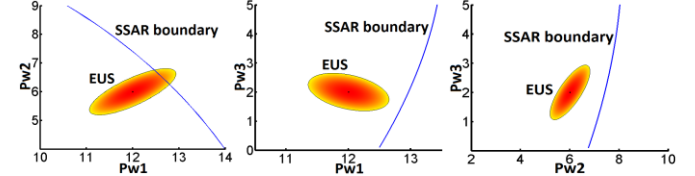

Fig. 8. Dimension-reduced SSAR boundary and EUS on different 2-D subspaces : (a) $(p_{w1},p_{w2})$; (b) $(p_{w1},p_{w3})$; (c) $(p_{w2},p_{w3})$

To justify the speculation above, a test is performed by reducing the $p_{wi0}$ ($i$=1,2,3) in four different directions, while keeping the amount of reduction constant as 0.2 p.u.. Tab. III lists the corresponding $P_{instab}$ in the original 3-D space after adjustments. The test results well verify the speculation since the $P_{instab}$ of Scenario #3 is larger than that of Scenario #1 but smaller than that of Scenario #2, while the $P_{instab}$ of Scenario #4 is larger than that before adjustment. The result reveals that wind generation curtailment does not necessarily facilitate the small-signal stability of power systems, which is opposite to common sense. In such a situation, the operators have to be choose the direction of curtailment carefully. In this sense, the proposed SSAR concept enables a graphic approach to support the operators' decision making on wind generation adjustment in a visual manner. However, to determine an optimal direction for wind generation adjustment in bulk power system operation, the current graphic tool is not enough, and some programming tools should be employed. In such a circumstance, the SSAR boundary could provide the programming with explicit constraints, remarkably reducing the computation complexity.

TABLE III
COMPARISON OF $P_{instab}$ AFTER CURTAILMENTS OF $p_{wi0}$ ($i$=1,2,3)

|  |  | $P_{instab}$ |
|---|---|---|
|  | Before adjustment | 6.66% |
| **Scenario #1** | Reduce $p_{w10}$ by 0.2 p.u | 3.36% |
| **Scenario #2** | Reduce $p_{w20}$ by 0.2 p.u | 4.34% |
| **Scenario #3** | Reduce $p_{w10}$, $p_{w20}$ by 0.1 p.u | 3.78% |
| **Scenario #4** | Reduce $p_{w30}$ by 0.2 p.u | 7.55% |

## VII. CONCLUDING REMARKS

This paper develops a conceptual tool, the SSAR, to geometrically depict the limit to accommodate uncertain WPIs without loss of small-signal stability of the bulk power system. We derive it based on the conventional SSSR concept and elucidate the maps among different nodal injection spaces. Theoretic analysis and simulations illustrate that the SSAR is capable of evaluating influence of uncertainty on small-signal stability from a region-wise point of view, leading to an intuitive but essential understanding on this challenging issue. This enables a simple way to measure and visualize stability boundaries in the uncertainty space of interest, which is highly desired by operators and engineers. Although the concept is derived here in terms of the issue of wind generation integration, it is quite general and can be directly applied to handle other uncertain renewable resources.

Some interesting research directions are open, including the fundamental theory, computational algorithms and applications. For the fundamental theory, it is natural to associate the SSAR with the power system flexibility, extending the traditional flexibility research from steady state security to system stability, and providing more insights on how to cope with uncertainty in operation. On the other hand, as mentioned in Part C, Section III, it is crucial albeit to develop efficient approaches to explicitly describe the boundaries of SSAR accurately in an enough large range. Another important research is to investigate the applications of SSAR to security monitoring, precaution and dispatch. For example, the robust small-signal stability boundary based on SSAR can be deployed as constraints in security-constrained unit commitment, economic dispatch, or optimal wind generation curtailment as mentioned in Section V.

ACKNOWLEDGEMENT

The authors would like to thank Felix F. Wu and Yunhe Hou for the helpful discussions.

APPENDIX

A. *3th-order model of synchronous generator with the simplified 3rd-order excitation model:*

$$\dot{\delta} = \omega_s(\omega - 1) \quad (A1)$$

$$T_J \dot{\omega} = -D\omega + P_m - P_e \quad (A2)$$

$$T'_{d0} \dot{e}'_q = -e'_q - (x_d - x'_d)i_d + e_{fq} \quad (A3)$$

$$T_B \dot{e}_{fq} = -e_{fq} + (T_A - T_C)v_R / T_A$$
$$\quad - K_A T_C v_M / T_A + K_A T_C (v_{ref} + v_S) / T_A \quad (A4)$$

$$T_A \dot{v}_R = K_A(v_{ref} + v_S - v_M) - v_R \quad (A5)$$

$$T_R \dot{v}_M = v_C - v_M \quad (A6)$$

B. *Beta distribution*

The probability density function of beta distribution is

$$f(x, \alpha, \beta, a, b) = \frac{1}{(b-a)B(\alpha, \beta)} \left(\frac{x-a}{b-a}\right)^{\alpha-1} \left(1 - \frac{x-a}{b-a}\right)^{\beta-1} \quad (B1)$$

where $a$ and $b$ are lower and upper limit of variable $x$. $B(\alpha,\beta)$ is the beta function:

$$B(\alpha, \beta) = \int_0^1 z^{\alpha-1}(1-z)^{\beta-1} dz \quad (B2)$$

where $z = (x-a)/(b-a)$. The mean $\mu$ and variance $\sigma^2$ of beta distribution are related to $\alpha$ and $\beta$ as follows:

$$\mu = \alpha(b-a)/(\alpha+\beta), \sigma^2 = \alpha\beta(b-a)^2 / \left((\alpha+\beta)^2(\alpha+\beta+1)\right) \quad (B3)$$